\documentclass[a4paper,12pt]{article}
\usepackage[english,russian]{babel}
\usepackage[T2A]{fontenc}
\usepackage[cp1251]{inputenc}
\usepackage[english]{babel}
\usepackage{amsthm}
\usepackage[tbtags]{amsmath}
\usepackage{amsfonts,amssymb}
\begin{document}

\renewcommand{\refname}{References}
\newtheorem{theorem}{Theorem}
\newtheorem{lemma}{Lemma}
\newtheorem{proposition}{Proposition}
\newtheorem{Cor}{Corollary}

\begin{center}
{\large\bf Centrally Essential Endomorphism Rings of Abelian Groups}
\end{center}
\begin{center}
Lyubimtsev O.V.\footnote{Nizhny Novgorod State University; email: oleg lyubimcev@mail.ru .}, Tuganbaev A.A.\footnote{National Research University `MPEI', Lomonosov Moscow State University; email: tuganbaev@gmail.com .}
\end{center}

\textbf{Abstract.}
We study Abelian groups $A$ with centrally essential endomorphism ring $\text{End}\,A$. If $A$ is a such group which is either a torsion group or a non-reduced group, then the ring $\text{End}\,A$ is commutative. We give examples of Abelian torsion-free groups of finite rank with non-commutative centrally essential endomorphism rings.

\textbf{MSC2010 database 16R99; 20K30}

\section{Introduction}
All rings considered are associative rings with non-zero identity element. A ring $R$ is said to be \textsf{centrally essential} if for any its non-zero element $a$, there exist two non-zero central elements $x,y\in R$ with $ax=y$. \footnote{It is clear that the ring $R$ with center $C$ is centrally essential if and only if the module $R_{C}$ is an essential extension of the module $C_{C}$.} Centrally essential rings are studied, for example, in \cite{MT18}, \cite{MT19}, \cite{MT19b}, \cite{MT19c}, \cite{MT20}, \cite{MT20b}, \cite{MT20c}.

It is clear that any commutative ring is centrally essential. In addition, every centrally essential semiprime ring is commutative; see \cite[Proposition 3.3]{MT18}. Therefore, in the study of centrally essential rings, we are only interested in non-commutative non-semiprime centrally essential rings.

Examples of non-commutative group algebras over fields are given in \cite{MT18}. 
For example, if $Q_8$ is the quaternion group of order $8$, then its group algebra over the field of order $2$ is a non-commutative centrally essential finite local ring of order $256$. 
In addition, in \cite{MT19}, it is proved that the Grassmann algebra of three-dimensional vector space over the field of order $3$ is a finite non-commutative centrally essential ring, as well. In \cite{MT19c}, there is an example of a centrally essential ring whose factor ring with respect to its prime radical is not a PI ring.

In Theorem 1.2(3) of this paper, we give an example of an Abelian torsion-free group of finite rank with centrally essential non-commutative endomorphism ring. In Example 3.9, we give additional examples of non-commutative centrally essential endomorphism rings of some Abelian torsion-free groups of infinite rank.

\textbf{1.1. Remark.} Let $A$ be an Abelian group which is either torsion group, or non-reduced group and let the endomorphism ring $\text{End}\,A$ be centrally essential. In Section 2 of this paper, we prove that the ring $\text{End}\,A$ is commutative. Therefore, when studying Abelian groups with non-commutative centrally essential endomorphism rings, only reduced torsion-free groups and reduced mixed groups are of interest.

Let $A$ be an Abelian torsion-free group with endomorphism ring $\text{End}\,A$ and let $\mathbb{Q}\text{End}\,A = \mathbb{Q}\otimes \text{End}\,A$ be the \textsf{quasi-endomorphism ring}\footnote{See 3.1 below.} of the group $A$. If the group $A$ has no a non-trivial quasi-decompositions\footnote{see 3.1 below.} then it is called \textsf{strongly indecomposable}. The \textsf{pseudo-socle} $\text{PSoc}\,A$ of the group $A$ is the pure subgroup of the group $A$ generated by all its minimal pure fully invariant subgroups.

The main result of this paper is Theorem 1.2.

\textbf{1.2. Theorem.} Let $A$ be a strongly indecomposable torsion-free Abelian group of finite rank, $\mathbb{Q}\text{End}\,A$ the quasi-endomorphism ring, and $A\neq \text{PSoc}\,A$. 

\textbf{1.} If $\mathbb{Q}\text{End}\,A$ is a centrally essential ring, then the ring $\mathbb{Q}\text{End}\,A/J(\mathbb{Q}\text{End}\,A)$ is commutative and $C(\mathbb{Q}\text{End}\,A) \cap M\neq 0$ for every minimal right ideal $M$ of $\mathbb{Q}\text{End}\,A$.

\textbf{2.} If the ring $\mathbb{Q}\text{End}\,A/J(\mathbb{Q}\text{End}\,A)$ is commutative, $\text{Soc}\,(\mathbb{Q}\text{End}\,A_{\mathbb{Q}\text{End}\,A}) = \text{Soc}\,(\mathbb{Q}\text{End}\,A_{C(\mathbb{Q}\text{End}\,A)})$ and 
$C(\mathbb{Q}\text{End}\,A) \cap M\neq 0$ for every minimal right ideal $M$ of $\mathbb{Q}\text{End}\,A$, then $\mathbb{Q}\text{End}\,A$ is a centrally essential ring.

\textbf{3.} Let $n > 1$ be an odd positive integer.  There exists a strongly indecomposable Abelian torsion-free group $A(n)$ of rank $2^n$ such that its endomorphism ring is a non-commutative centrally essential ring.
 
For convenience, we give some definitions and notation used in the paper. The necessary ring-theoretical information not listed in the paper can be found in \cite{Row88}. The necessary information on Abelian groups not listed in the paper can be found in \cite{Fuc15} and \cite{KMT03}.

If $R$ is a ring, then we denote by $C(R)$ and $J(R)$ the center and the Jacobson radical of the ring $R$, respectively. We use the additive form for Abelian groups. We denote by $\text{End}\,A$ the endomorphism ring of the Abelian group $A$. If 
$A = \bigoplus_{p\in P} A_p$ is a decomposition of the torsion Abelian group $A$ into the direct sum of $p$-components, then  $\text{supp} \,A = \{p\in P \mid A_p\neq 0\}$. We use the following notation: $\mathbb {Z}_{p^k}$ (resp., $Z_{p^k}$) is the residue ring (resp., the additive group) modulo $p^k$;~ $\mathbb Q$ (resp., $Q$) is the ring (resp., the additive group) of rational numbers; $Z_{p^{\infty}}$ is a quasi-cyclic Abelian group; $\hat{\mathbb Z}_p$ is the ring of $p$-adic integers. If $A$ is an Abelian torsion-free group, then $\mathbb{Q}\text{End}\,A$ and $\text{PSoc}\,A$ are the quasi-endomorphism ring and the pseudo-socle of the group $A$, respectively.

A ring $R$ is said to be \textsf{local} if the factor ring $R/J(R)$ is a division ring.

For a module $M$, the \textsf{socle} $\text{Soc}\,M$ is the sum of all simple submodules of $M$; if $M$ does not contain a simple submodule, then $\text{Soc}\,M = 0$ by definition. 

An Abelian group $A$ is said to be \textsf{divisible} if $nA = A$ for any positive integer $n$. An Abelian group is said to be \textsf{reduced} if it does not contain a non-zero divisible subgroup and \textsf{non-reduced}, otherwise. 

A subgroup $B$ of the Abelian group $A$ is said to be \textsf{pure} if the equation $nx = b\in B$, which has a solution in the group $A$, has a solution in $B$.

\section{Non-Reduced Abelian Groups with Centrally Essential Endomorphism Rings} 

\textbf{2.1. Lemma.} Let $A$ be a module and $A = \bigoplus_{i\in I}A_i$ a direct decomposition of the module $A$. The endomorphism ring $\text{End} A$ is centrally essential if and only if for every $i\in I$ the following conditions hold.

\textbf{1)} $A_i$ is a fully invariant submodule in $A$; 

\textbf{2)} the ring $\text{End} A_i$ is centrally essential. 

\textbf{Proof.} Let $\text{End}\,A = E$ be a centrally essential ring. If condition 1) does not hold and $A_i$ is not a fully invariant submodule for some $i\in I$, then there exists a subscript $j\in I$, $j\neq i$, such that $\text{Hom}\,(A_i, A_j) = e_jEe_i\neq 0$, where $e_i$ and $e_j$ are the projections from the module $A$ onto the submodules $A_i$ and $A_j$, respectively. In addition,
$$
e_i\cdot e_jEe_i = 0\neq e_jEe_i = e_jEe_i\cdot e_i,
$$
i.e., the idempotent $e_i$ is not central; this contradicts to \cite[Lemma 2.3]{MT18}. 

If every $A_i$ is a fully invariant submodule in $A$, $i\in I$, then $\text{End}\,A\cong \text{End}\,A_i\times \text{End}\,\overline{A}_i$, where $\overline{A}_i$ is a complement direct summand of $A_i$. It is obvious that if $\text{End}\,A_i$ is not centrally essential ring, then and  $\text{End}\,A$ is not a centrally essential ring.

If conditions 1) and 2) hold, then $\text{End}\,A\cong \prod_{i\in I}\text{End}\,A_i$ and each of the ring $\text{End}\,A_i$ is centrally essential. It is clear that the ring $\text{End}\,A$ is centrally essential, as well.~\hfill$\square$

\textbf{2.2. Lemma.} The endomorphism ring of a divisible Abelian group $A$ is centrally essential if and only if either $A\cong Q$ or $A\cong Z_{p^{\infty}}$.

\textbf{Proof.} Let $A = F(A)\bigoplus T(A)$, where $0\neq F(A)$ is the torsion-free part and $0\neq T(A)$ is the torsion part of the group $A$. Then $F(A)$ is not a fully invariant subgroup in $A$ (see \cite[Theorem 7.2.3]{Fuc15}) and, by Lemma 2.1, the ring $\text{End}\,A$ is not centrally essential. Hypothetically $F(A)$ or $T(A)$ is a direct sum of ${\mathbb{Z}}_p^{\infty}$ or $\mathbb{Q}$. Clearly, if the number of terms is $>1$, the ring $\text{End}\,A$ has a noncentral idempotent which gives a contradiction.~\hfill$\square$
 
Let $A = \bigoplus_{p\in P} A_p$ be the decomposition of the torsion Abelian group $A$ into the direct sum of its primary components. It follows from Lemma 2.1 that
$\text{End}\,A$ is a centrally essential ring if and only if each of the ring $\text{End}\,A_p$ is centrally essential.

\textbf{2.3. Lemma.} The endomorphism ring of a primary Abelian group $A_p$ is centrally essential if and only if $A_p\cong Z_{p^k}$ or $A_p\cong Z_{p^{\infty}}$.

\textbf{Proof.} If the group $A_p$ is not indecomposable, then it has a co-cyclic direct summand; see \cite[Corollary 5.2.3]{Fuc15}. By considering \cite[Theorem 7.1.7, Example 7.1.3 and Theorem 7.2.3]{Fuc15}, this summand or summands complement to it are not fully invariant in $A$. Consequently, $A_p\cong Z_{p^k}$ or $A_p\cong Z_{p^{\infty}}$. The converse is obvious, since rings $\mathbb{Z}_{p^k}$ and $\hat{\mathbb{Z}}_p$ are commutative.~\hfill$\square$

\textbf{2.4. Theorem.} Let $A = D(A)\bigoplus R(A)$ be a non-reduced Abelian group, where $0\neq D(A)$ and $0\neq R(A)$ are the divisible part and the reduced part of the group $A$, respectively. The endomorphism ring of the group $A$ is centrally essential if and only if 
$A = D(A)\bigoplus R(A)$, where $R(A) = \bigoplus_{p\in P'}Z_{p^k}$ and $D(A)\cong Q$ or $D(A)\cong \bigoplus_{p\in P''} Z_{p^{\infty}}$;
$P', P''$ are the subsets of different primes with $P' \cap P'' = \emptyset$.

\textbf{Proof.} Let $\text{End}\,A$ be a centrally essential ring. We verify that $D(A)$ and $R(A)$ are fully invariant subgroups in $A$. Indeed, it is well known that $\text{Hom}\,(D(A), R(A)) = 0$. Next, if $R(A)$ is a torsion-free group, then $\text{Hom}\,(R(A), D(A)) \neq 0$ (see \cite[Theorem 7.2.3]{Fuc15}); this contradicts to Lemma 2.1. It is also clear that $\text{Hom}\,(R(A), D(A)) \neq 0$ if $R(A)$, $D(A)$ are torsion groups and $\text{supp} \,R(A) \cap \text{supp} \,D(A) \neq \emptyset$. It follows from Lemma 2.3 that $R(A)$ is the direct sum of its cyclic $p$-components and it follows from Lemma 2.2 that $D(A)\cong Q$ or $D(A)\cong \bigoplus_{p\in P} Z_{p^{\infty}}$.

The converse assertion directly follows from Lemmas 2.1, 2.2 and 2.3.~\hfill$\square$

\textbf{2.5. Corollary.} The endomorphism ring of a non-reduced Abelian group is centrally essential if and only if the ring is commutative.

\textbf{Proof.} Indeed, it follows from Theorem 2.4 that any centrally essential endomorphism ring of a non-reduced Abelian group is the direct product of rings whose components can be only the rings $\mathbb{Z}_{p^k}$, $\mathbb{Q}$ and $\hat{\mathbb{Z}}_p$. ~\hfill$\square$

It follows from Corollary 2.5 that only reduced Abelian groups can have non-commutative centrally essential endomorphism rings.

\section{Proof of Theorem 1.2} 

\textbf{3.1. Quasi-decompositions and strongly indecomposable groups.}\\ Let $A$ and $B$ be two Abelian torsion-free groups. One says that $A$ is \textsf{quasi-contained} in $B$ if $nA\subseteq B$ for some positive integer $n$. If $A$ is quasi-contained in $B$ and $B$ is quasi-contained in $A$ (i.e., if $nA\subseteq B$ and $mB\subseteq A$ for some
$n,m\in \mathbb{N}$), then one says that $A$ is \textsf{quasi-equal} to $B$ (we write $A\doteq B$). A quasi-equality $A\doteq\bigoplus_{i\in I}A_i$ is called a \textsf{quasi-decomposition} (or a \textsf{quasi-direct decomposition}) of the Abelian group $A$; these subgroups $A_i$ are called \textsf{quasi-summands} of the group $A$. If the group $A$ does not have non-trivial quasi-decompositions, then $A$ is said to be \textsf{strongly indecomposable}. A ring $\mathbb{Q}\otimes \text{End}\,A$ is called the \textsf{quasi-endomorphism ring} of the group $A$; we denote it by $\mathbb{Q}\text{End}\,A$; see details in \cite[Chapter I, \S 5]{KMT03}. We note that
$$
\mathbb{Q}\text{End}\,A = \{\alpha\in \text{End}_{\mathbb{Q}}(Q\otimes A) \mid (\exists n\in \mathbb{N}) (n\alpha\in \text{End}\,A)\}.
$$
It is well known (e.g., see \cite [Proposition 5.2]{KMT03}) that the correspondence
$$
A \doteq e_1A\bigoplus \ldots \bigoplus e_kA \to \mathbb{Q}\text{End}\,A = \mathbb{Q}\text{End}\,A e_1\bigoplus \ldots \bigoplus\mathbb{Q}\text{End}\,A e_k
$$
between finite quasi-decompositions of the torsion-free group $A$ and finite decompositions of the ring $\mathbb{Q}\text{End}\,A$ in to a direct sum of left ideals, where
$\{e_i \mid i = 1,\ldots,k\}$ is a complete orthogonal system of idempotents of the ring $\mathbb{Q}\text{End}\,A$, is one-to-one.

\textbf{3.2. Proposition.} \footnote{cf. \cite[Proposition 2.2]{MT20c}}
The endomorphism ring $E$ of an Abelian torsion-free group $A$ is centrally essential if and only if the quasi-endomorphism ring $\mathbb{Q}E$ of $A$ is centrally essential.

\textbf{Proof.} Let $0\neq \tilde{a}\in \mathbb{Q}E$. For some $n\in \mathbb{N}$, we have $n\tilde{a} = a\in E$ and there exist  $x,y\in C(E)$ with $ax = y\neq 0$. In this case, $\tilde{a}\tilde{x} = \tilde{y}$, where  $\tilde{x} = x$, $\tilde{y} = \cfrac{1}{n}\cdot y\in C(\mathbb{Q}E)$, i.e., $\mathbb{Q}E$ is a centrally essential ring.

Conversely, for every $0\neq a\in E$, there exist non-zero $\tilde{x},\tilde{y}\in C(\mathbb{Q}E)$ with $a\tilde{x} = \tilde{y}$. In addition, there exist $n, m\in \mathbb{N}$ such that $n\tilde{x}\in C(E)$ and $m\tilde{y}\in C(E)$. Then $ax = y$, where $x = mn\tilde{x}, y = mn\tilde{y}\in C(E)$.~\hfill$\square$

Let $A\doteq \bigoplus_{i = 1}^n A_i = A'$ be a decomposition of the Abelian torsion-free group $A$ of finite rank into a quasi-direct sum of strongly indecomposable groups (e.g., see \cite[Theorem 5.5]{KMT03}). By considering Lemma 2.1 and Proposition 3.2, we obtain that the ring $\text{End}\,A$ is centrally essential if and only if all subgroups $A_i$ are fully invariant in $A'$, and every ring $\text{End}\,A_i$ is centrally essential. 
Therefore, the problem of describing Abelian torsion-free groups of finite rank with centrally essential endomorphism rings is reduced to the similar problem for strongly indecomposable groups.

\textbf{3.3. Proposition.} Let $A$ be a strongly indecomposable Abelian group and $A = \text{PSoc}\,A$. The ring $\text{End}\,A$ is centrally essential if and only if 
$\text{End}\,A$ is a commutative ring.

\textbf{Proof.} If $A = \text{PSoc}\,A$, then $\text{End}\,A$ is a semiprime ring (e.g., see \cite [Theorem 5.11]{KMT03}). It follows from \cite[Proposition 3.3]{MT18} that the ring $\text{End}\,A$ is commutative. The converse is obvious.~\hfill$\square$

\textbf{3.4. Proposition.} Let $R$ be a local Artinian ring which is not a division ring and $C(R) = C$. 

\textbf{1.} If $R$ is a centrally essential ring, then the ring $R/J(R)$ is commutative and
$C \cap M\neq 0$ for every minimal right ideal $M$.

\textbf{2.} If the ring $R/J(R)$ is commutative, $\text{Soc}\,(R_C) = \text{Soc}\,(R_R)$ and $C \cap M\neq 0$ for every minimal right ideal $M$, then $R$ is a centrally essential ring.

\textbf{Proof.} \textbf{1.} It is well known that if $R$ is an Artinian ring, then $J(R)$ is a nilpotent ideal of some index $k$. We note that if $M$ is a minimal right ideal of $R$, then $MJ(R) = 0$. Indeed, if $MJ(R) = M$, then
$$
M = MJ(R) = MJ^2(R) = \ldots = MJ^k(R) = 0.
$$
Let $R$ be a centrally essential ring. It follows from \cite [Theorem 2]{MT19} that the ring $R/J(R)$ is commutative. 

We assume that $C \cap M = 0$ for some minimal right ideal $M$. By assumption, for $0\neq a\in M$, there exist $x, y\in C(R)$ with $ax = y\neq 0$. Since $x\notin J(R)$ (otherwise, $ax = 0$), we have that $x$ is invertible in $R$ and $a = x^{-1}y\in C$, a contradiction. 

\textbf{2.} Let $C \cap M\neq 0$ for every minimal right ideal $M$ of $R$. We verify that $M\subseteq C$. Let $C\cap M = K$. By assumption, the 
ring $R/J(R)$ is commutative; therefore, we have $rs - sr\in J(R)$, for all $r, s\in R$. For every $k\in K$, we have $k(rs - sr) = 0$. On the other hand, since $k\in C$, we have that $(kr)s = ksr = s(kr)$ and $kr\in C$. In addition, $kr\in M$. Consequently, $K$ is a right ideal. From the property that $M$ is a minimal right ideal, we have that $K = M$ or $K = 0$. However $K\neq 0$; therefore, $K = M$ and $M\subset C$. Therefore, $\text{Soc}\,(R_C) = \text{Soc}\,(R_R)\subseteq C$. It follows from \cite [Theorem 3]{MT19} that $R$ is a centrally essential ring. ~\hfill$\square$

\textbf{3.5. Example.} We will find centrally essential endomorphism rings of strongly indecomposable Abelian torsion-free groups of rank 2 and 3.

If $A$ is an strongly indecomposable group of rank $2$, then the ring $\text{End}\,A$ is commutative (e.g., see \cite[Theorem 4.4.2]{Fat07}). Consequently, $\text{End}\,A$ is a centrally essential ring. Let $A$ be a strongly indecomposable group of rank $3$. Then the algebra $\mathbb{Q}\text{End}\,A$ is isomorphic to one of the following $\mathbb{Q}$-algebras (\cite[Theorem 2]{Ch98}):
$$
K\cong\left\{ \left(\begin{matrix}
x & 0 & z\\0 & x & 0\\0 & 0 & x 
\end{matrix}\right) \mid x, z\in \mathbb{Q}\right\},\quad R\cong\left\{ \left(\begin{matrix}
x & y & z\\0 & x & 0\\0 & 0 & x 
\end{matrix}\right) \mid x, y, z\in \mathbb{Q}\right\},
$$
$$
S\cong\left\{ \left(\begin{matrix}
x & y & z\\0 & x & ky\\0 & 0 & x 
\end{matrix}\right) \mid x, y, z\in \mathbb{Q},\,0\neq k\in \mathbb{Q},\,k = const\right\},
$$
$$
T\cong\left\{ \left(\begin{matrix}
x & y & z\\0 & x & t\\0 & 0 & x 
\end{matrix}\right) \mid x, y, z, t\in \mathbb{Q}\right\}.
$$
The rings $K$, $R$, $S$ are commutative; consequently, they are centrally essential. The ring $T$ is not commutative (in addition, $\text{PSoc}\,A$ is of rank $1$). We have 
$$
J(T) = \left\{ \left(\begin{matrix}
0 & y & z\\0 & 0 & t\\0 & 0 & 0 
\end{matrix}\right) \mid y, z, t\in \mathbb{Q}\right\},
$$
$$
C(T) = \left\{ \left(\begin{matrix}
x & 0 & z\\0 & x & 0\\0 & 0 & x 
\end{matrix}\right) \mid x, z \in \mathbb{Q}\right\},
$$ 
$$
M = \left\{ \left(\begin{matrix}
0 & 0 & 0\\0 & 0 & t\\0 & 0 & 0 
\end{matrix}\right) \mid t\in \mathbb{Q}\right\},
$$ 
where $M$ is the minimal right ideal of $T$. We note that the ring $T/J(T)$ is commutative, but $C(T) \cap M = 0$. 
It follows from Proposition 3.4(1) that the ring $T$ is not centrally essential.
As a result, we obtain that endomorphism rings strongly indecomposable groups of rank 2 or 3 are centrally essential if and only if they are commutative.

\textbf{3.6. Example.} Let $V$ be a vector $\mathbb{Q}$-space with basis $e_1, e_2, e_3$ and let $\Lambda(V)$ be the Grassmann algebra of the space $V$, i.e., $\Lambda(V)$ is an algebra with operation $\wedge$, generators $e_1, e_2, e_3$ and defining relations
$$
e_i \wedge e_j + e_j \wedge e_i = 0 \quad \mbox {for of all} \quad i, j = 1, 2, 3.
$$
Then $\Lambda(V)$ is a $\mathbb{Q}$-algebra of dimension 8 with basis $\{1, e_1, e_2, e_3, e_1 \wedge e_2, e_2 \wedge e_3, e_1 \wedge e_3, e_1 \wedge e_2 \wedge e_3\}$ and $\Lambda(V)$ is a non-commutative centrally essential ring (see details in \cite[Example 1]{MT19}). 
We consider the regular representation of the algebra $\Lambda(V)$. If $x\in \Lambda(V)$,
$$
x = q_0\cdot 1 + q_1e_1 + q_2e_2 + q_3e_3 + q_4e_1 \wedge e_2 + q_5e_2 \wedge e_3 + q_6e_1 \wedge e_3 + q_7e_1 \wedge e_2 \wedge e_3,
$$
then the matrix $A_x\in Mat_8(\mathbb{Q})$ has the form
$$
\left(\begin{matrix}
q_0 & q_1 & q_2 & q_3 & q_4 & q_5 & q_6 & q_7\\
0 & q_0 & 0 & 0 & -q_2 & 0 & -q_3 & q_5\\
0 & 0 & q_0 & 0 & q_1 & -q_3 & 0 & -q_6\\
0 & 0 & 0 & q_0 & 0 & q_2 & q_1 & q_4\\
0 & 0 & 0 & 0 & q_0 & 0 & 0 & q_3\\
0 & 0 & 0 & 0 & 0 & q_0 & 0 & q_1\\
0 & 0 & 0 & 0 & 0 & 0 & q_0 & -q_2\\
0 & 0 & 0 & 0 & 0 & 0 & 0 & q_0 \end{matrix}\right).
$$
We denote by $R$ the corresponding subalgebra in $\text{Mat}_8(\mathbb{Q})$. It is clear that the radical $J(R)$ consists of properly upper triangular matrices in $R$ and $A_x\in C(R)$ if and only if $q_1 = q_2 = q_3 = 0$. In addition, 
$Soc (R_R) = \{A_x = (a_{ij})\in R \mid a_{ij} = 0, i\neq 1, j\neq 8\}$ and
$Soc (R_C) = \{A_x = (a_{ij})\in C(R) \mid a_{ii} = 0\}$.
Since $Soc (R_C)\neq Soc (R_R)$, the corresponding condition of Proposition 3.4(2) is not necessary. Therefore, we obtain the negative answer to \cite[Open questions 4.5(2)]{MT19}.

\textbf{3.7. The completion of the proof of Theorem 1.2.}\\
\textbf{1, 2.} It is known that the ring $\mathbb{Q}\text{End}\,A$ is a local Artinian ring (e.g., see \cite[Corollary 5.3]{KMT03}). It remains to use Proposition 3.4.

\textbf{3.} In \cite[Proposition 2.5.]{MT19}, it is proved that the Grassmann algebra $\Lambda(V)$ over a field $F$ of characteristic $0$ or $p\neq 2$ is a centrally essential ring if and only if the dimension of the space $V$ is odd. We set $F = \mathbb{Q}$. It is known (e.g., see \cite{PV83}) that every $\mathbb{Q}$-algebra of dimension $n$ can be realized as the quasi-endomorphism ring of an Abelian torsion-free group of rank $n$. Therefore, by considering Example 3.6 and Proposition 3.2, we obtain the required property.~\hfill$\square$
 
Under conditions of Theorem 1.2, if the rank of the group $A$ is square-free, then the ring $\mathbb{Q}\text{End}\,A/J(\mathbb{Q}\text{End}\,A)$ is commutative \cite[Lemma 4.2.1]{Fat07}. By considering Proposition 3.2, we obtain

\textbf{3.8. Corollary.} Let $A$ be a strongly indecomposable Abelian torsion-free group of finite rank, $A\neq \text{PSoc}\,A$ and the rank of the group $A$ is square-free. 

\textbf{1.} If the endomorphism ring $\text{End}\,A$ of the group $A$ is centrally essential, then $C(\mathbb{Q}\text{End}\,A) \cap M\neq 0$ for every minimal right ideal $M$ of $\mathbb{Q}\text{End}\,A$.

\textbf{2.} If $\text{Soc}\,(\mathbb{Q}\text{End}\,A_{\mathbb{Q}\text{End}\,A}) = \text{Soc}\,(\mathbb{Q}\text{End}\,A_{C(\mathbb{Q}\text{End}\,A)})$ and $C(\mathbb{Q}\text{End}\,A) \cap M\neq 0$ for every minimal right ideal $M$ of $\mathbb{Q}\text{End}\,A$, then $\text{End}\,A$ is a centrally essential ring.

\textbf{3.9. Example.} Let $R = \mathbb{Z}[x,y]$ be the polynomial ring in two variables $x$ and $y$. We use the construction described in \cite[Proposition 7]{J16}. We consider the ring
$$
T(R) = \left\{ \left(\begin{matrix}
f & d_1(f) & g\\0 & f & d_2(f)\\0 & 0 & f 
\end{matrix}\right) \mid f, g\in \mathbb{Z}[x,y]\right\},
$$
where $d_1$, $d_2$ are two derivations of the ring $\mathbb{Z}[x,y]$, $d_1 = \cfrac{\partial}{\partial x}$, $d_2 = \cfrac{\partial}{\partial y}$.
Then $T(R)$ is a non-commutative ring with $J(R) = e_{13}R\subseteq C(T(R))$, where $e_{13}$ is the matrix unit; see \cite[Corollary 8]{J16}. If $0\neq a\in T(R)\setminus C(T(R))$, then $0\neq ae_{13}\in C(T(R))$. Therefore, $T(R)$ is a centrally essential ring.
Since $T(R)$ is a countable ring with reduced torsion-free additive group, it follows from the familiar Corner theorem (e.g., see \cite[Theorem 29.2]{KMT03}) that for the ring $T(R)$, there exist $\mathfrak{M}$ of Abelian groups $A_i$ such that $\text{End}\,A_i\cong T(R)$ 
and $\text{Hom}\,(A_i, A_j) = 0$ for all $i\neq j$, where $\mathfrak{M}$ is an arbitrary preset cardinal number; see \cite{DG82}, \cite{CG85}. We note that the endomorphism ring of the direct sum of such groups is a non-commutative centrally essential ring, as well.

\section{Remarks and Open Questions}

\textbf{4.1. Open question.} Is it true that there exist strongly indecomposable Abelian torsion-free groups of rank $< 8$ whose endomorphism rings are non-commutative centrally essential rings?

\textbf{4.2. Open question.} An Abelian group is said to be \textsf{super-decomposable} if it does not have non-zero indecomposable direct summands. Is it true that there exist a super-decomposable Abelian group with non-commutative centrally essential endomorphism ring?

\textbf{4.3. Open question.} Is it true that the endomorphism ring of the direct sum of all the groups $A(n)$ from Theorem 1.2(3) is a non-commutative centrally essential ring with polynomial identity?

\textbf{4.4. Open question.} Is it true that there exists an Abelian group $A$ with centrally essential endomorphism ring $\text{End}\,A$ which is not a ring with polynomial identity?

\textbf{4.5.} A ring is said to be \textsf{right distributive} (resp., \textsf{right uniserial}) if the lattice its right ideals is distributive (resp., is a chain). If the endomorphism ring $\text{End}\,A$ of the group $A$ of finite rank is right uniserial, then $\text{End}\,A$ is an invariant\footnote{A ring is said to be \textsf{invariant} if all its one-sided ideals are ideals.} principal right ideal domain; see \cite[Proposition 3.4]{Fat88}. Therefore, if $A$ is an Abelian torsion-free group of finite rank and $\text{End}\,A$ is a centrally essential right uniserial ring, then $\text{End}\,A$ is commutative (\cite[Proposition 2.8]{MT19b}). In addition, it is known that every right distributive local ring is right uniserial. Consequently, every centrally essential right distributive quasi-endomorphism ring of an Abelian torsion-free group of finite rank is commutative. We note that there exist non-commutative uniserial Artinian centrally essential rings; see \cite{MT20}.

In connection to 4.5, we formulate open questions 4.6 and 4.7.

\textbf{4.6. Open question.} Is it true that there exist Abelian torsion-free groups of finite rank whose endomorphism rings are non-commutative right distributive centrally essential rings? 

\textbf{4.7. Open question.} Is it true that there exist Abelian groups whose endomorphism rings are non-commutative right distributive (or right uniserial) centrally essential rings?

\end{document}